\documentclass[11pt, a4paper]{article}

\usepackage{amsmath,amsfonts}
\usepackage{xcolor}
\usepackage{url}
\usepackage{amsthm}
\usepackage[hang]{footmisc}
\usepackage{authblk}
\usepackage{comment}
\usepackage{mathrsfs}
\usepackage[pagebackref]{hyperref}
\usepackage{bbm}
\usepackage{enumerate}
\numberwithin{equation}{section}
\newtheorem{theorem}{Theorem}[section]
\newtheorem{lemma}[theorem]{Lemma}
\newtheorem{proposition}[theorem]{Proposition}

\newtheorem{definition}[theorem]{Definition}
\theoremstyle{remark}
\newtheorem{remark}[theorem]{Remark}
\newtheorem{observation}[theorem]{Observation}

\newcommand{\W}{W_{\xi,\varepsilon}}
\newenvironment{pf}{{\noindent\bf Proof.}}{\qed\newline}

\newcommand{\R}{\mathbb{R}}
\usepackage{etoolbox}
\makeatletter
\newcounter{author}
\renewcommand*\author[1]{%
  \stepcounter{author}%
  \ifnum\c@author=1
    \gdef\@author{#1}%
  \else
    \xdef\@author{\unexpanded\expandafter{\@author\and#1}}%
  \fi
  \csgdef{author@\the\c@author}{#1}}
\newcommand*\email[1]{%
  \csgdef{email@\the\c@author}{#1}}
\newcommand*\address[1]{%
  \csgdef{address@\the\c@author}{#1}}
\AtEndDocument{%
  \xdef\author@count{\the\c@author}%
  \c@author=1
  \print@authors}
\newcommand*\print@authors{%
  \ifnum\c@author>\author@count
  \else
    \print@author{\the\c@author}%
    \advance\c@author by 1
    \expandafter\print@authors
  \fi}
\newcommand*\print@author[1]{%
  \par\medskip
  \begin{tabular}{@{}l@{}}%
    \textsc{Address of \csuse{author@#1}}\\
    \csuse{address@#1}\\
    \textit{E-mail address}:
    \href{mailto:\csuse{email@#1}}{\csuse{email@#1}}
  \end{tabular}}
  \newcommand{\subjclass}[2][1991]{%
  \let\@oldtitle\@title%
  \gdef\@title{\@oldtitle\footnotetext{#1 \emph{Mathematics subject classification.} #2}}%
}
\makeatother
\title{A multiplicity result for a double perturbed Schr\"odinger-Bopp-Podolsky-Proca system}
\date{\vspace{-5ex}}
\author{Matteo Talluri}
\address{Dipartimento di Matematica\\Università di Bologna\\ Piazza di Porta S. Donato, 5\\ 40126 Bologna, Italy.}
\email{matteo.talluri@unibo.it}

\usepackage{lipsum}
\usepackage{titlefoot}
\begin{document}
\maketitle
\unmarkedfntext{Date: November 4, 2023\newline2010 \emph{Mathematics Subject Classification}. 35J40, 35Q55, 53C21.\newline \emph{Keywords: Nonlinear Schr\"odinger equation, Bopp-Podolsky-Proca electromagnetic theory, Calculus of variations, Riemannian manifolds.}}
\begin{abstract}\noindent In this paper we will prove a multiplicity result for a double perturbed Schr\"odinger-Bopp-Podolsky-Proca system on a compact 3-dimensional Riemannian manifold without boundary. We will prove, using the Lusternik-Schnirelmann Category, that the number of solutions of the system depends on the topological properties of the manifold.\end{abstract}
\section{Introduction.}
We look for positive solutions of the Schr\"odinger-Bopp-Podolsky-Proca system \begin{equation}\label{sbpp}\left\{\begin{array}{l}
-\varepsilon^2 \Delta_g u+ u+q^2 \phi u=|u|^{p-2} u \\
-\varepsilon^2\Delta_g \phi+ \varepsilon^4\Delta_g^2 \phi+\phi=4 \pi u^2
\end{array} \quad \text { in } M\right.\end{equation}
where $M$ is a 3-dimensional compact Riemannian manifold without boundary, $q$ is a positive parameter, $p\in(4,6)$ and $\Delta_g$ is the Laplace-Beltrami operator. 
The above system can be derived considering the motion of a charged particle $\psi$ in an electromagnetic field $(\mathbf{A},\phi)$ generated by itself in the Bopp-Podolsky-Proca setting. The Bopp-Podolsky model, unlike Maxwell's one, solves the problem of the finiteness of the electron's energy. While in the Maxwell model, a charged particle generates an electrostatic field that is proportional to $|x|^{-1}$, in the Bopp-Podolsky model it is proportional to $
\frac{1-e^{-\frac{1}{|x|}}}{|x|},$ and therefore its gradient has finite $L^2$ norm. A further generalisation of the theory, that is due to Proca, is that the electromagnetic field is assumed to have mass $m$. The motion of the particle $\psi$ is described by the Schr\"odinger equation, whose lagrangian is \begin{equation}\label{lagrangiana schrodinger}
\mathcal{L}_S(\psi)=i\hbar\partial_t\psi\overline{\psi}-\frac{\hbar^2}{2m^2}|\nabla_g\psi|^2+\frac{2}{p}|\psi|^p.
\end{equation}If we want to consider the effect of the electromagnetic field on such particle we can use the so-called coupling rule, which consists in replacing $\partial_t$ and $\nabla_g$ in (\ref{lagrangiana schrodinger}) with \[
\partial _t\rightarrow\partial_t+i\frac{q}{\hbar}\phi\quad,\quad\nabla_g\rightarrow\nabla_g-i\frac{q}{\hbar c}\mathbf{A}.
\]Therefore the new lagrangian of $\psi$ is \[
\mathcal{L}_S(\psi,\mathbf{A},\phi)=i\hbar\partial_t\psi\overline{\psi}-q\phi|\psi|^2-\frac{\hbar^2}{2m^2}|\nabla_g\psi-i\frac{q}{\hbar c}\mathbf{A}|^2+\frac{2}{p}|\psi|^p
.\]The complete lagrangian $\mathcal{L}_{SBBP}$ of the system is obtained by adding to $\mathcal{L}_S$ the contribute of the electromagnetic field $\mathcal{L}_{BPP}$ in the Bopp-Podolsky theory and with the assumption, that is due to Proca, that the electromagnetic field has mass $m$. We refer to \cite{hebey2019electro} for the complete derivation of the system. If we take the Euler-Lagrange equations of $\mathcal{L}_{SBBP}$ in the case when the magnetic field vanishes (that is $\mathbf{A}=0$) and we look for a stationary solution of the form $\psi(x,t)=u(x)e^{i\omega t}$ where $u\ge 0$ we get, after a rescaling and normalizing some parameters, system (\ref{sbpp}), see \cite{hebey2020strong} for the details. From a mathematical point of view, the first to study the Bopp-Podolsky model were D'Avenia and Siciliano in \cite{MR3957980}, where an existence result and some properties of the solutions are proved. These kinds of problems have been widely studied by Hebey on compact manifolds (see \cite{MR4396239} and \cite{hebey2019electro}) and by a wide range of authors in $\R^3$. We mention, without claiming to be exhaustive, some of them: \cite{MR4652600}, \cite{MR4648641}, \cite{MR4636514}, \cite{MR4547444}, and \cite{MR4519835}. 
We want to show that if $\varepsilon$ is small enough, there exist at least as many weak solutions with low energy of the system (\ref{sbpp}) as the Lusternik–Schnirelmann category of the manifold; whose Definition we recall to be: \begin{definition}
Let $X$ be a topological space and consider a closed subset $A \subset X$. We say that $A$ has category $k$ relative to $X$ ($\operatorname{cat}_X A=k$ ) if $A$ is covered by $k$ closed sets $A_j, j=1, \ldots, k$, which are contractible in $X$, and $k$ is the minimum integer with this property. We simply denote $\operatorname{cat} X=\operatorname{cat}_X X$.
\end{definition}Moreover, we will show that at least another solution with a higher energy exists. In particular, the main result of this paper is the following Theorem.\begin{theorem}\label{risultato principale}
For all $p\in(4,6)$ and $\varepsilon$ small enough, there exist at least $\operatorname{cat} M+1$ non-constant positive solution of the system (\ref{sbpp}).
\end{theorem}
The proof of the Theorem is divided into two parts. The first one relies on the \emph{Photography methods}, which is a topological tool that was firstly introduced by Bahari and Coron \cite{bahri1987nonlinear} for a non-linear elliptic equation on bounded domains and then re-adapted by Benci and Cerami \cite{benci1991effect} and Benci, Cerami and, Passaseo \cite{MR1205376} both for systems and on compact manifolds. The main idea of this method is that for every point of the manifold we can construct a solution of the system with low energy concentrated near this point. Conversely, we can show that every function whose energy functional is small enough is concentrated near a point of the manifold. Thanks to these results we can show that the category of the set of solutions with low energy is bounded from below by the category of the manifold. Since the manifold is compact this result also shows that there exists at least two solutions of the system. The second part of the proof relies on a method that was used by Benci, Bonanno and Micheletti in \cite{benci2007multiplicity}; the main idea is to construct a contractible set of positive functions.  This will imply that the energy functional admits a critical point that is not a minimum. This method has been applied to study a wide variety of systems similar to (\ref{sbpp}), for instance, Ghimenti and Micheletti in \cite{MR3392985} used it to study a Schr\"odinger- Maxwell type system while Siciliano, Figueiredo, and Mascaro in \cite{figueiredo2023multiple} and \cite{mascaro2022positive} used a variant of the Photography method to investigate the solutions of a Schr\"odinger-Bopp-Podolsky equation (that is our system where $m=0$). It is important to point out that system (\ref{sbpp}) has already been studied by D'Avenia and Ghimenti in \cite{d2022multiple} when $\varepsilon$ appears only in the first equation. The fact that the perturbation appears only in the first equation means that the limit problem they obtain consists in only one equation whose solution is unique and well-known, which is not true in our case. The lack of uniqueness will force us to prove that almost all the energy of one of such solutions is concentrated in a ball of radius big enough. The paper is organised into seven sections. Section \ref{sezione 2} introduces notations and preliminary definitions, section \ref{sezione 3} contains some remarks about the Nehari manifold, and section \ref{sezione 4} is dedicated to the study of the limit problem. Sections \ref{sezione 5} and \ref{sezione 6} are dedicated to the proof of the first part of Theorem \ref{risultato principale}, in particular, the connection between the topology of the manifold and low energy functions is exploited. Finally, in section \ref{sezione 7} we conclude the proof Theorem \ref{risultato principale} showing that there exists at least another solution with higher energy.
\section{Preliminaries.}\label{sezione 2}
In the rest of the paper, for every $\varepsilon>0$ fixed, we will use the following rescaled norms, which are equivalent to the standard ones of $H^2(M)$, $H^1(M)$ and $L^p(M)$, respectively:$$
\|v\|_{2,\varepsilon}^2:=\varepsilon\int_M|\Delta_gv|^2d\mu_g+\frac{1}{\varepsilon} \int_M\left|\nabla_g v\right|^2 d \mu_g+\frac{1}{\varepsilon^3} \int_M v^2 d \mu_g,
$$ $$
\|v\|_{1,\varepsilon}^2:=\frac{1}{\varepsilon} \int_M|\nabla_g v|^2 d \mu_g+\frac{1}{\varepsilon^3} \int_M v^2 d \mu_g, \quad|v|_{p, \varepsilon}^p:=\frac{1}{\varepsilon^3} \int_M|v|^p d \mu_g.
$$By a standard computation, it is easy to see that there exists $C$ that does not depend on $\varepsilon$, such that $$
|v|_{p, \varepsilon} \leq C\|v\|_{1,\varepsilon}
\quad \forall p\in[2,6].$$The following result, which can be found in \cite{hebey2019electro}, will allow us to reduce the study of the system to that of a single equation.
\begin{lemma}\label{proprietà seconda equazione}
For every $u\in H^1(M)$ there exists a unique $\phi(u)\in H^4(M)\cap C^2(M)$ solution of \begin{equation}
-\varepsilon^2\Delta_g v+\varepsilon^4 \Delta_g^2 v+v=4 \pi u^2 \quad \text { in } M.
\end{equation}Moreover\begin{enumerate}[(i)]
    \item there exists $C>0$ independent on $\varepsilon$, such that for every $u\in H^1(M)$, $ \|\phi(u)\|_{H^2}\le C|u|_2^2$ and $\|\phi(u)\|_{H^4}\le C|u|_4^2;$
\item if $\varepsilon$ is small enough then $\phi(u)\ge 0;$
\item for every $t\neq0$ we have $\phi(tu)=t^2\phi(u).$
\end{enumerate}
\end{lemma}
By this result we can write (\ref{sbpp}) as \begin{equation}\label{singola equazione}
-\varepsilon^2 \Delta_g u+u+\phi(u) u=|u|^{p-2} u \quad \text { in } M
.\end{equation}At this point we can introduce a functional whose Euler-Lagrange equation is (\ref{singola equazione}). To do this, we will use the following Lemma:\begin{lemma}\cite[Lemma 2.3]{d2022multiple}\label{derivata termine misto}
    For every $u\in H^1(M)$ we define the functional \[
    G(u)=\int_M u^2 \phi(u) d \mu_g
    ,\]then $G$ is $C^1$ and we have \[
    G'(u)[h]=4\int_{M}\phi(u)uhd\mu_g\quad\forall\:h\in H^1(M).
    \]Moreover, if $u_n\to \bar{u}$ in $H^1(M)$ then $G(u_n)\to G(\bar{u})$.
\end{lemma}Thanks to Lemma \ref{derivata termine misto} it is easy to see that positive solutions of (\ref{singola equazione}) are critical points of the the functional \[
J_\varepsilon(u)=\frac{1}{2}\|u\|_{1,\varepsilon}^2+\frac{1}{4\varepsilon^3}\int_M\phi(u)u^2\:d\mu_g-\frac{1}{p}|u^+|_{p,\varepsilon}^p
\quad u\in H^1(M).\] We conclude this section with some remarks about the manifold $M$. Thanks to Nash's embedding Theorem \cite{nash1956imbedding}, we can assume that $M$ is isometrically embedded in $\R^n$ for some $n>3.$ Moreover, since $M$ is compact, there exist a radius $r>0$ such for every $\xi\in M $ the exponential map \[\operatorname{exp}_{\xi}:B(0,r)\subseteq T_pM\to B_g(\xi,r)\subseteq M
\]is a diffeomorphism. In addition, see for instance \cite{lang2012differential}, for every $y\in B(0,r)$, the following expansions of the metric $g$ hold: \begin{equation}
    \label{espansioni tensore metrico} \left(g_{\xi}\right)_{i j}(y)=\delta_{i j}+O(|y|)\quad\left|g_{\xi}(y)\right|:=\operatorname{det}(\left(g_{\xi})_{i j}\right)(y)=1+O(|y|).
\end{equation}
\section{Remarks on Nehari manifold.}\label{sezione 3}
Since positive solutions of (\ref{singola equazione}) are critical points of the functional $J_\varepsilon$ we can consider the Nehari manifold \begin{equation}\label{nehari}\
\mathcal{N}_\varepsilon:=\{u\in H^1(M)\setminus
\{0\}\::\:N_\varepsilon(u)=0\}
\end{equation}where \begin{equation}\label{nehari funzionale}
N_\varepsilon(u):=J'_\varepsilon(u)[u]=\|u\|_{1,\varepsilon}^2+\frac{1}{\varepsilon^3}\int_M\phi(u)u^2\:d\mu_g-|u^+|_{p,\varepsilon}^p.\end{equation}
\begin{remark}
If $u$ belongs to the set $\mathcal{N}\varepsilon$, we can observe the following relation:
\begin{equation*}
N_{\varepsilon}(u) = J_{\varepsilon}^{\prime}(u)[u] = \|u\|_{1,\varepsilon}^2 + \frac{1}{\varepsilon^3} \int_M \phi(u) u^2d\mu_g - \left|u^{+}\right|_{p, \varepsilon}^p = 0.
\end{equation*}
Thus, we obtain the equation:
\begin{equation}\label{remark 1}
\|u\|_{1,\varepsilon}^2 = -\frac{1}{\varepsilon^3} \int_M \phi(u) u^2 d\mu_g + \left|u^{+}\right|_{p, \varepsilon}^p.
\end{equation}
By substituting this expression into $J_\varepsilon$, we find:
\begin{equation}\label{J=Gamma su N}
J_\varepsilon(u) = \left(\frac{1}{2} - \frac{1}{p}\right) \frac{1}{\varepsilon^3} \int_M \left|u^{+}\right|^p d\mu_g- \frac{1}{4 \varepsilon^3} \int_M u^2 \phi(u)d\mu_g = \int_ M \Gamma(u)d\mu_g,
\end{equation}
where
\begin{equation}\label{definizione gamma}
\Gamma(u) := \left(\frac{1}{2} - \frac{1}{p}\right) \frac{1}{\varepsilon^3} \left|u^{+}\right|^p - \frac{1}{4 \varepsilon^3} u^2 \phi(u).
\end{equation}
Similarly, from equation (\ref{remark 1}), we have:
\begin{equation}\label{N uguale a 0}
|u^+|_{p,\varepsilon}^p = \|u\|_{1,\varepsilon}^2 + \frac{1}{\varepsilon^3} \int_M \phi(u)u^2  d\mu_g,
\end{equation}
which, combined with the expression for $J_\varepsilon$, leads to:
\begin{equation}\label{espressione J senza norma p}
J_\varepsilon(u) = \left(\frac{1}{2} - \frac{1}{p}\right) \|u\|_{1,\varepsilon}^2 + \left(\frac{1}{4} - \frac{1}{p}\right) \int_M \phi(u)u^2 d\mu_g.
\end{equation}
Furthermore, we can also find:
\begin{equation}
\label{J senza misto}J_\varepsilon(u) = \frac{\|u\|_{1,\varepsilon}^2}{4} + \left(\frac{1}{4} - \frac{1}{p}\right) |u^+|_{p,\varepsilon}^p.
\end{equation}
\end{remark}
\begin{lemma}\cite[Lemma 3.1]{d2022multiple}\label{proprietà tempo}
 If $p>4$, then:\begin{enumerate}[(i)]
     \item for every $u \in H^1(M)$ with $u^{+} \not \equiv 0$, there exists a unique $t_u>0$ such that $t_u u \in\mathcal{N}_\varepsilon$;
\item the map $u \in H^1(M)  \mapsto t_u \in(0,+\infty)$ is continuous;
\item there exists $C>0$ such that, for every $\varepsilon>0$ and $u \in \mathcal{N}_{\varepsilon}$, $\left|u^{+}\right|_{p, \varepsilon} \geq C$;
\item for every $u \in \mathcal{N}_{\varepsilon}$, $ N_{\varepsilon}^{\prime}(u) \neq 0$;
\item there exists $C>0$ such that, for every $\varepsilon>0$, $m_{\varepsilon}:=\inf _{\mathcal{N}_{\varepsilon}} J_{\varepsilon} \geq C>0$. \end{enumerate}
\end{lemma}
\begin{lemma}\cite[Lemma 3.2]{d2022multiple}\label{palais smale}
    Let $\varepsilon>0.$ Then, for every $c>0,$ the functional $J_\varepsilon$ satisfies the Palais-Smale condition at level $c$. Moreover, if $\{u_n\}$ is a Palais-Smale sequence for $J_\varepsilon$ restricted to $\mathcal{N}_\varepsilon$, then it is also a Palais-Smale sequence for $J_\varepsilon.$
\end{lemma}
The proof of these two results relies on classical techniques and can be found in detail in the cited papers, so we prefer to omit it here.
\section{The limit problem.}\label{sezione 4}Consider the following problem in the whole space:\begin{equation}\label{limit problem}
\begin{cases}-\Delta u+u+ u \phi=|u|^{p-2} u & \text { in } \mathbb{R}^3 \\ -\Delta \phi+\Delta^2\phi+\phi=4\pi u^2 & \text { in } \mathbb{R}^3 \\ u>0 & \text { in } \mathbb{R}^3 \end{cases}
.\end{equation} and as before define $\phi(u)$ as the solution of the second equation. If we consider the functional \begin{equation}
\label{funzionale limite}
J_{\infty}(u):=\frac{1}{2}\|u\|_{H^1}^2+\frac{1}{4} \int_{\R^3}u^2\phi(u)-\frac{1}{p}\left|u^{+}\right|_p^p,\quad u\in H^1(\R^3)
\end{equation}
and the Nehari manifold $$
\mathcal{N}_{\infty}:=\left\{u \in H^1\left(\mathbb{R}^3\right) \backslash \{0\}: J_{\infty}^{\prime}(u)[u]=0\right\}
,$$ it is possible to prove (see for instance \cite{MR2843921}) that \[
m_\infty:=\inf_\mathcal{N_\infty}J_\infty
\]is attained, however, the uniqueness of a such minimizer is not known. 
In the rest of the paper, $U$ will be a fixed solution of the first equation in (\ref{limit problem}) that realizes the minimum. We will call $U_\varepsilon(\cdot)=U(\cdot/\varepsilon)$ and, for all $\xi\in M$ we define \begin{equation}\label{soluzione in carta}
W_{\xi, \varepsilon}(\cdot):=U_{\varepsilon}\left(\exp _{\xi}^{-1} \cdot\right) \chi_r\left(\left|\exp _{\xi}^{-1} \cdot\right|\right)
.\end{equation}\begin{lemma}\label{proprietà soluzione limite}
    For every $\xi\in M$ the following facts hold:\begin{enumerate}[(i)]
        \item $\displaystyle\lim _{\varepsilon \rightarrow 0}\left\|W_{\xi, \varepsilon}\right\|_{1,\varepsilon}^2=\|U\|_{H^1}^2$;
        \item $\displaystyle\lim _{\varepsilon \rightarrow 0}\left|W_{\xi, \varepsilon}\right|_{q, \varepsilon}^q=|U|_q^q$ for $q \in[1,6] ;$
        \item $\displaystyle\lim _{\varepsilon \rightarrow 0} \frac{1}{\varepsilon^3} \int_M W_{\xi, \varepsilon}^2 \phi(W_{\xi, \varepsilon})\:d\mu_g=\int_{\R^3}U^2\phi(U);$
        \item $\displaystyle\lim _{\varepsilon \rightarrow 0} t_{W_{\xi, \varepsilon}}=1$.
    \end{enumerate}Where $U$ is a minimizer of (\ref{funzionale limite}) on $\mathcal{N}_\infty.$
\end{lemma}
\begin{pf}
For the proof of $(i)$ and $(ii)$ we refer to \cite[Lemma 4.1]{d2022multiple}. In order too prove $(iii)$ we define \[
\widetilde\phi_\varepsilon(z)=\phi(W_{\xi,\varepsilon})(\text{exp}(\varepsilon z))\chi_r(|\varepsilon z|)
.\] Using the expansions of the metric tensor (\ref{espansioni tensore metrico}) it is easy to see that 
$\|\widetilde\phi_\varepsilon\|_{H^2(\R^3)}^2\le C\|\phi(W_{\xi,\varepsilon})\|_{2,\varepsilon}^2,$ and by $(i)$ in Lemma \ref{proprietà seconda equazione} we get that \[
\|\phi(W_{\xi,\varepsilon})\|_{2,\varepsilon}^2\le C|W_{\xi,\varepsilon}|_\varepsilon^2\le C\|U\|_{H^1}^2,
\]so $\widetilde \phi_\varepsilon$ is bounded in $H^2$ and therefore it converges to a function $\widetilde\phi$ weakly in $H^2(\R^3)$ and strongly in $H^1(\R^3).$ We want to show that, in a weak sense, \[
-\Delta \widetilde\phi+\Delta^2\widetilde\phi+\widetilde\phi=4\pi U^2.
\]Given $f\in C^\infty_c(\R^3)$, we have that $\text{sup}f\subseteq B(0,\frac{r}{\varepsilon})$ if $\varepsilon$ is small enough. Therefore, if we define \[
f_\varepsilon(x)=f\bigg(\frac{\text{exp}_\xi^{-1}(x)}{\varepsilon}\bigg),\]we have that $f_\varepsilon\in C^{\infty}_c(M)$ and $\text{sup}f_\varepsilon\subseteq B_g(\xi,r)$. Since $\phi(W_{\xi,\varepsilon})$ solves the second equation in (\ref{sbpp}) we have that \[
\int_M\frac{\nabla_g\phi(W_{\xi,\varepsilon})\nabla_g f_\varepsilon}{\varepsilon}+\varepsilon\Delta_g\phi(W_{\xi,\varepsilon})\Delta_g f_\varepsilon+\frac{\phi(W_{\xi,\varepsilon})f_\varepsilon}{\varepsilon^3}
\:d\mu_g=\frac{4\pi}{\varepsilon^3}\int_MW_{\xi,\varepsilon}^2f_\varepsilon\:d\mu_g.\]If we pass in coordinates and make the change of variables $y=\varepsilon z$ we get that the LHS is equal to \begin{align*}
    &\int_{B(0,\frac{r}{\varepsilon})}g^{ij}(\varepsilon z)\partial_i\widetilde\phi_\varepsilon(z)\partial_jf(z)|g_{i,j}(\varepsilon z)|^{1/2}\:dz+\int_{B(0,\frac{r}{\varepsilon})}\widetilde\phi_\varepsilon(z)f(z)|g_{i,j}(\varepsilon z)|^{1/2}\:dz\\+&\int_{B(0,\frac{r}{\varepsilon})}\frac{\partial_i(|g_{i,j}(\varepsilon z)|^{1/2}g^{i,j}(\varepsilon z)\partial_j\widetilde\phi_\varepsilon(z)\partial_i(|g_{i,j}(\varepsilon z)|^{1/2}g^{i,j}(\varepsilon z)\partial_jf(z))}{|g_{i,j}(\varepsilon z)|^{1/2}}\:dz
\end{align*}
and therefore, by (\ref{espansioni tensore metrico}), when $\varepsilon\to0$ it converges to \[
\int_{\R^3}\nabla\widetilde\phi\nabla f+\widetilde\phi f-\Delta\widetilde\phi\Delta f\:dz.\]In the same way we get\[\frac{4\pi}{\varepsilon^3}
\int_MW^2_{\xi,\varepsilon}f_\varepsilon\:d\mu_g=4\pi\int_{B(0,\frac{r}{\varepsilon})}
|U(z)|^2\chi_r(|\varepsilon z|)f(z)|g_{i,j}(\varepsilon z)|^{1/2}\:dz\to4\pi\int_{\R^3}|U(z)|^2f(z)\:dz,\]and so $\widetilde \phi=\phi(U).$ To conclude it is enough to observe that \[
\frac{1}{\varepsilon^3}\int_MW^2_{\xi,\varepsilon}\phi(W_{\xi,\varepsilon})\:d\mu_g=\int_{\R^3}U^2(z)\chi(|\varepsilon z|)^2\widetilde\phi_\varepsilon(z)|g_{i,j}(\varepsilon z)|\:dz\to\int_{\R^3}|U(z)|^2\phi(U)\:dz.
\]To prove $(iv)$, we start observing that $t^2_{W_{\xi,\varepsilon}}$ is the unique solution of\[
0=t^2\|W_{\xi,\varepsilon}\|_{1,\varepsilon}^2+\frac{t^4}{\varepsilon^3}\int_M\phi(W_{\xi,\varepsilon})W_{\xi,\varepsilon}^2\:d\mu_g-t^p|W_{\xi,\varepsilon}^+|^p_{p,\varepsilon}
\]that can be written as \[
0=\|W_{\xi,\varepsilon}\|_{1,\varepsilon}^2+\frac{t_{W_{\xi\varepsilon}}^2}{\varepsilon^3}\int_M\phi(W_{\xi,\varepsilon})W_{\xi,\varepsilon}^2\:d\mu_g-t_{W_{\xi,\varepsilon}}^{p-2}|W_{\xi,\varepsilon}^+|^p_{p,\varepsilon}.
\]Since $p>4$ and by $(i)$, $(ii)$ and $(iii)$ we get that $t_{W_{\xi,\varepsilon}}$ is bounded and therefore it converges to some ${t}$ which solves \[
\|U\|_{H^1}^2=-{t}^2\int_{\R^3}U^2\phi(U)+{{t}}^{p-2}|U|_p^p
.\]But this equation has a unique solution in $t$, and since $U$ solves (\ref{limit problem}), it holds that \[
\|U\|_{H^1}^2=-\int_{\R^3}U^2\phi(U)+|U|_p^p
,\]and so ${t}=1.$ \qedhere
\end{pf} 
\section{The map \texorpdfstring{$\Psi_\varepsilon$}{Ψε}.}\label{sezione 5}
If we choose an $U$ that minimizes (\ref{funzionale limite}) on $\mathcal{N}_\infty$, for every $\varepsilon>0$ we can define a map $$\Psi_\varepsilon:M\to\mathcal{N}_\varepsilon$$ $$\xi\to t_{W_{\xi,\varepsilon}}\W,$$ since $t_u$ depends continuously on $u\in H^1(M)$ and $U$ is fixed we have that the above map is continuous. Using equation (\ref{espressione J senza norma p}) and Lemma \ref{proprietà soluzione limite} the following Propositions easily follow. The details of the proofs can be found in \cite[Proposition 4.3 and Lemma 5.1]{MR3392985}.
\begin{proposition}\label{N limitato su psi}For every $\delta>0$ there exists an $\varepsilon_0=\varepsilon_0(\delta)$ such that, if $\varepsilon<\varepsilon_0$, then $J_\varepsilon(\Psi_\varepsilon(\xi))<m_\infty+\delta.$
\end{proposition}
\begin{proposition}
    There exists $\alpha>0$ such that, for every $\varepsilon>0$ it holds that $m_\varepsilon\ge\alpha.$\label{minimo lontano da 0}
\end{proposition}
\begin{observation}
    \label{limsup}By Proposition \ref{N limitato su psi} it also follows that $$
\limsup _{\varepsilon \rightarrow 0} m_\varepsilon\leq m_{\infty}.
$$
\end{observation}
\section{The barycenter map.}\label{sezione 6}For every $u\in \mathcal{N}_\varepsilon$ we can define a point $\beta(u)\in\R^n$ by $$
\beta(u)=\frac{\int_M x \Gamma(u) d \mu_g}{\int_M \Gamma(u) d \mu_g}
,$$ where $\Gamma(u)$ was defined in $(\ref{definizione gamma})$. We have to see that $\beta$ is well defined and to do so is enough to observe that if $u\in \mathcal{N}_\varepsilon$, we can use equation (\ref{J=Gamma su N}) to infer that $J_\varepsilon(u)=\int_M\Gamma (u)$ and by Proposition \ref{minimo lontano da 0} we have $\Gamma(u)\ge m_\varepsilon>0.$
In the following, we will need to consider a \emph{good} partition of our manifold $M$, whose Definition is:
\begin{definition}
Given $\varepsilon>0$ we will say that a finite partition $P_\varepsilon=\{P_j^\varepsilon\}$, ${j\in \Lambda_\varepsilon}$ of $M$ is good if:\begin{enumerate}[(i)]
    \item $P_j^\varepsilon$ is closed for every $j\in \Lambda_\varepsilon;$
    \item $P_j^\varepsilon\cap P_l^\varepsilon\subset\partial P_j^\varepsilon\cap\partial P_j^\varepsilon$ for all $j\neq l;$
    \item There exists $C>0$ and $r_1(\varepsilon)\ge r_2(\varepsilon)>C\varepsilon$  such that there exists $q_j^\varepsilon\in P_j^\varepsilon$ for which \[
    B_g(q_j^\varepsilon,\varepsilon)\subset P_j^\varepsilon\subseteq B_g(q_j^\varepsilon,r_2(\varepsilon))\subset B_g(q_j^\varepsilon,r_1(\varepsilon))
    ;\]
    \item There exists $\nu(M)\in \mathbb{N}$ such that every $q\in M$ is contained in at most $\nu(M)$ balls $B_g(q_j^\varepsilon,r_1(\varepsilon)).$
\end{enumerate}
\end{definition}
The existence of such partition was proved in \cite{benci2007multiplicity}. Moreover, the following result, which was originally proved in \cite[Lemma 5.3]{benci2007multiplicity} and then adapted to the Bopp-Podolsky setting in \cite[Lemma 4.4]{d2022multiple}, prevents the vanishing of $u$ on the Nehari manifold.
\begin{proposition}\label{stima partizione}
    There exists $\gamma>0$ such that, for every $\delta>0$ and $\varepsilon\in(0,\varepsilon_0(\delta))$ where $\varepsilon_0(\delta)$ is the one of Proposition \ref{N limitato su psi}, if $P_\varepsilon$ is a good partition of $M$ then for every $u\in\mathcal{N}_\varepsilon$ there exists $\Bar{j}\in\Lambda_\varepsilon$ such that \[
    \frac{1}{\varepsilon^3}\int_{P_{\bar{j}}^\varepsilon}|u^+|^p\:d\mu_g\ge\gamma.
    \]
\end{proposition}
\begin{proposition}\label{q(u)}
    There exists $\eta\in(0,1)$ and $\delta_0\in(0,m_\infty)$ such that for all $\delta\in(0,\delta_0)$, $\varepsilon\in(0,\varepsilon_0(\delta))$ and $u\in\mathcal{N}_\varepsilon\cap J_\varepsilon^{m_\infty+\delta}$ there exists $q:=q(u)\in M$ such that \[
    \int_{B_g(q,\frac{r}{2})}\Gamma(u)\:d\mu_g\ge(1-\eta)m_\infty,
    \]where $\varepsilon(\delta)$ is as in Proposition \ref{N limitato su psi}.
\end{proposition}
\begin{pf} First of all we prove the statement for $u\in\mathcal{N}_\varepsilon\cap J_\varepsilon^{m_\varepsilon+\delta}.$ Assume by contradiction that it is false. Therefore there exists $\delta_k,\varepsilon_k$ and $u_k\in\mathcal{N}_{\varepsilon_k}$ such that: \begin{equation}
\lim_{k\to+\infty}\delta_k=\lim_{k\to+\infty}\varepsilon_k=0
\end{equation}\begin{equation}
    m_{\varepsilon_k}\le J_{\varepsilon_k}(u_k)\le m_{\varepsilon_k}+2\delta_k
\end{equation}\begin{equation}
    \int_{B_g(q,\frac{r}{2})}\Gamma(u)d\mu_g\le(1-\eta)m_\infty, 
\end{equation}moreover, by Observation \ref{limsup} we can also assume that \[
m_{\varepsilon_k}+2\delta_k\leq m_\infty+3\delta_k\le 2m_\infty
.\]Since $\delta_k\to0$, $u_k\in\mathcal{N}_{\varepsilon_k}$ and $\phi(u_k)\ge0$ we have that \begin{equation}\label{bound H1}
2m_\infty\ge\left(\frac{1}{2}-\frac{1}{p}\right)\|u_k\|_{1,\varepsilon_k}^2+\left(\frac{1}{4}-\frac{1}{p}\right)\frac{1}{\varepsilon_k^3}\int_M\phi(u_k)u_k^2\:d\mu_g\ge\left(\frac{1}{2}-\frac{1}{p}\right)\|u_k\|_{1,\varepsilon_k}^2.
\end{equation}By the Ekeland principle we can assume that there exists $\sigma_k\to 0$ such that \begin{equation}\label{ekland}
|J'_{\varepsilon_k}(u_k)(\phi)|\le \sigma_k\|\phi\|_{1,\varepsilon_k},
\end{equation}and by Proposition \ref{stima partizione} there exists $P_k^{\varepsilon_k}\in P_{\varepsilon_k}$ such that \begin{equation}\label{w positiva}
\frac{1}{\varepsilon_k^3}\int_{P_k^{\varepsilon_k}}|u_k^+|^p\:d\mu_g\ge\gamma>0.
\end{equation}If $q_k$ belongs to a such $P_k^{\varepsilon_k}$ we define \[
w_k(z)=\begin{cases}u_k(\exp_{q_k}(\varepsilon_kz))\chi(\varepsilon_k|z|)&\quad \text{if $z\in B(0,\frac{r}{\varepsilon_k})$}\\0&\quad \text{if $z\notin B(0,\frac{r}{\varepsilon_k})$}
\end{cases},\]it holds that $w_k\in H^1_0(\R^3)$ and by equation (\ref{bound H1}) \[
\left\|w_k\right\|_{H^1\left(\mathbb{R}^3\right)}^2 \leq C\left\|u_k\right\|_{1,\varepsilon_k}^2 \leq C .
\]So there exists $w\in H^1_0(\R^3)$ such that $w_k\rightharpoonup w$ in $H^1_0(\R^3)$ and $w_k\to w$ in $L^t_{\text{loc}}(\R^3)$ for every $2\le t\le6.$ Let $\phi(u_k)$ be the solution of the second equation in (\ref{sbpp}), if we define \[
\Phi_k(z)=\phi(u_k)(\exp_{q_k}(\varepsilon_kz)\chi(\varepsilon_k|z|))
\]and we argue as in the third point of Lemma \ref{proprietà soluzione limite} we get that $\Phi_k\rightharpoonup\phi(w)$ in $H^2(\R^3)$ and $\Phi_k\to\phi(w)$ in $W^{1,p}_{\text{loc}}(\R^3)$ for all $2\le p\le6.$ Arguing as in \cite[Proposition 5.3]{MR3392985} and using (\ref{w positiva}) we can show that $w$ is a ground state solution of (\ref{limit problem}), and \begin{equation}
\liminf_{k\to+\infty}m_{\varepsilon_k}=m_\infty\label{liminf}.
\end{equation}
To conclude we choose $T>0$ such that \[
w_k\to w \quad\text{and}\quad \Phi_k\to\phi(w)\quad\text{in}\:L^p(B(0,T))
,\]$k$ such that $\varepsilon_kT\le\frac{r}{2}$, and \begin{equation}\label{energy}
\int_{B(0,T)}\left[\left(\frac{1}{2}-\frac{1}{p}\right)(w_k^+(z))^p-\frac{1}{4}w_k^2(z)\Phi_k(z)\right]g_{i,j}(\varepsilon z)\ge(1-\frac{\eta}{2})m_\infty.
\end{equation}We observe that:
\begin{align*}
    (1-\frac{\eta}{2})m_\infty&\le\int_{B(0,T)}\left[\left(\frac{1}{2}-\frac{1}{p}\right)(w_k^+(z))^p-\frac{1}{4}w_k^2(z)\Phi_k(z)\right]g_{i,j}(\varepsilon z)\\
    &=\int_{B_g(q_k,\varepsilon_kT)}\left(\frac{1}{2}-\frac{1}{p}\right)(u_k^+)^p-\frac{1}{4}u_k^2\phi(u_k)\:d\mu_g\\
    &=\int_{B_g(q_k,\varepsilon_kT)}\Gamma(u_k)d\mu_g\le\int_{B_g(q_k,\frac{r}{2})}\Gamma(u_k)\:d\mu_g\le(1-\eta)m_\infty
\end{align*}that is a contradiction. So we proved the statement for $u\in\mathcal{N}_\varepsilon\cap J_\varepsilon^{m_\varepsilon+2\delta}.$
By observation \ref{limsup} and equation (\ref{liminf}) we get \[
\lim_{\varepsilon\to0}m_\varepsilon=m_\infty,
\]hence, if $\varepsilon$ and $\delta$ are small enough $J_\varepsilon^{m_\infty+\delta}\subseteq J_\varepsilon^{m_\varepsilon+2\delta}$, and the statement follows.
\end{pf}
\begin{observation}
    The proof of Proposition \ref{q(u)} would have been easier if the limit problem (\ref{funzionale limite}) had a unique solution, which in our case could not be true. 
    However, even if we don't know $w$, the functions $w_k$ converge to some $w$ in $L^p_{\text{loc}}(\R^n)$. Also, $w$ is a ground state, therefore we can choose $T$ such that (\ref{energy}) holds. After this choice, we can take $k$ large enough such that $\varepsilon_kT\le \frac{r}{2}$ and thanks to the rescaling we conclude.
\end{observation}
\begin{proposition}\label{omotopa identità}
    There exists $\delta_0\in(0,m_\infty)$ such that for every $\delta\in(0,\delta_0)$, $\varepsilon\in(0,\varepsilon(\delta))$ and $u\in\mathcal{N}_\varepsilon\cap J_\varepsilon^{m_\infty+\delta}$ it holds that $\beta(u)\in M_r(m))$ where $\varepsilon(\delta)$ is defined in Proposition \ref{N limitato su psi} and \[
    M_{r(M)}=\{x\in \R^n\::\:d(x,M)<r(M)\}.
    \]Moreover $\beta\circ\Phi_\varepsilon:M\to M_{r(M)}$ is homotopic to the immersion.
\end{proposition}
\begin{pf}
    Given $u\in\mathcal{N}_\varepsilon\cap J_\varepsilon^{m_\infty+\delta}$ and $\eta\in(0,1)$, if $\varepsilon$ and $\delta$ are small enough, by Proposition \ref{q(u)} we can find a $q:=q(u)\in M$ such that \[
    \int_{B_g(q(u),\frac{r}{2})}\Gamma(u)\:d\mu_g \ge(1-\eta)m_\infty.
    \]Since $u\in\mathcal{N}_\varepsilon\cap J_\varepsilon^{m_\infty+\delta}$ it also holds that \[
    J_\varepsilon(u)=\int_M\Gamma(u)\:d\mu_g\le m_\infty+\delta
    \]and hence \begin{align*}
    |\beta(u)-q|&\le\frac{\left|\int_M(x-q)\Gamma(u)\:d\mu_g\right|}{\int_M\Gamma(u)\:d\mu_g}\\
    &\le\frac{\left|\int_{B_g(q,\frac{r}{2})}(x-q)\Gamma(u)\:d\mu_g\right|}{\int_M\Gamma(u)\:d\mu_g}+\frac{\left|\int_{M\setminus B_g(q,\frac{r}{2})}(x-q)\Gamma(u)\:d\mu_g\right|}{\int_M\Gamma(u)\:d\mu_g}\\&\le\frac{r}{2}+\text{diam}(M)\left(1-\frac{(1-\eta)m_\infty}{m_\infty+\delta}\right).
    \end{align*}Since \[
    \frac{(1-\eta)m_\infty}{m_\infty+\delta}\to 1 \quad\text{if}\quad \eta,\delta\to 0,
    \]we have \[
    |\beta(u)-q|\le r
    \]if $\eta $ and $\delta$ are small enough.
\end{pf}
We are finally ready to prove the first part of Theorem \ref{risultato principale}. To do this, we will use two well-known results.
\begin{theorem}\cite[Theorem 9.10, Lemma 7.10]{ambrosetti2007nonlinear}\label{categoria punti critici}
     Let $J$ be a $C^{1,1}$ real functional on a complete $C^{1,1}$ manifold $\mathcal{N}$. If $J$ is bounded from below and satisfies the Palais Smale condition, then it has at least $\operatorname{cat}\left(J^d\right)$ critical point in $J^d$ where $$J^d=\{u \in \mathcal{N}: J(u) \leq d\}.$$ Moreover if $\mathcal{N}$ is contractible and $\operatorname{cat} (J^d)>1$, there exists at least one critical point $u \notin J^d.$
\end{theorem}
\begin{proposition}\cite{benci1991effect}\label{categoria composizione}
   Let $X_1$ and $X_2$ be topological spaces. If $g_1: X_1 \rightarrow X_2$ and $g_2: X_2 \rightarrow X_1$ are continuous operators such that $g_2 \circ g_1$ is homotopic to the identity on $X_1$, then cat $X_1 \leq$ cat $X_2$.
\end{proposition}
To prove the first statement of Theorem \ref{risultato principale} is enough to use Lemma \ref{palais smale} to get that $J_\varepsilon$ fulfils the hypothesis of Theorem \ref{categoria punti critici} on $\mathcal{N}_\varepsilon\cap J_\varepsilon^{m_\infty+\delta}$, therefore $J_\varepsilon$ has at least $\operatorname{cat}(\mathcal{N}_\varepsilon\cap J_\varepsilon^{m_\infty+\delta})$ critical points. By Propositions \ref{q(u)} and \ref{categoria composizione} we also have $$\operatorname{cat}(\mathcal{N}_\varepsilon\cap J_\varepsilon^{m_\infty+\delta})\ge\operatorname{cat}(M)> 1,$$where the last inequality follows by the fact that $M$ is compact.
\section{High energy solutions.}\label{sezione 7}
In this section, we conclude the proof of Theorem \ref{risultato principale} proving that the system (\ref{sbpp}) has at least a solution with higher energy. Following \cite[Section 6]{benci2007multiplicity} is enough to show that there exists a set $T_\varepsilon\subseteq\mathcal{N}_\varepsilon$ with the following properties:\begin{enumerate}[(i)]
    \item $\Psi_\varepsilon(M)\subseteq T_\varepsilon$;
    \item Every $u\in T_\varepsilon$ is positive;
    \item $T_\varepsilon$ is contractible in $\mathcal{N}_\varepsilon\cap J_\varepsilon^{c_\varepsilon}$ where $c_\varepsilon$ is bounded by a suitable constant $c$ that does not depend on $\varepsilon$.
\end{enumerate}
Indeed, by the same argument of the proof of Theorem \ref{risultato principale}, we can show that \[\operatorname{cat}\left(T_\varepsilon\cap J_\varepsilon^{m_\infty+\delta}\right)>1,\]
moreover since $T_\varepsilon$ is contractible, by the second part of Theorem \ref{categoria punti critici} $J_\varepsilon$ admits a critical point $u\in T_\varepsilon$ such that $J_\varepsilon(u)> m_\infty+\delta.$ 
In order to construct the set we take $V\in H^{1}(\R^3)$ a positive function and for a fixed $\xi_0\in M$ we define \[
v_{\varepsilon}:=V_{\varepsilon}\left(\exp _{\xi_0}^{-1} \cdot\right) \chi_r\left(\left|\exp _{\xi_0}^{-1} \cdot\right|\right)\]and the cone \[
C_{\varepsilon}:=\big\{u=\theta v_{\varepsilon}+(1-\theta) W_{\xi, \varepsilon},\: \theta \in[0,1], \:\xi \in M\big\}
\]where $W_{\xi, \varepsilon}$ is defined in (\ref{soluzione in carta}). It is easy to see that $C_\varepsilon$ is contractible in $H^1(M)$; moreover, by Proposition \ref{proprietà soluzione limite}, it is also compact in $H^1(M).$ By Proposition \ref{proprietà tempo} the function $u\to t_u$ is continuous, hence the set \[
T_{\varepsilon}:=\left\{t_u u: u \in C_{\varepsilon}\right\}
\]is still compact and contractible. Clearly $T_\varepsilon\subseteq\mathcal{N}_\varepsilon$ and since $t_u$, $v_\varepsilon$ and $W_{\xi,\varepsilon}$ are positive, $C_\varepsilon$ contains only positive functions. If we define \[
c_\varepsilon:=\max _{u \in C_\varepsilon} J_\varepsilon\left(t_u u\right)=\max_{v\in T_\varepsilon}J_\varepsilon(v)
,\]we have that \[
T_\varepsilon\subseteq \mathcal{N}_\varepsilon \cap J_\varepsilon^{c_\varepsilon},
\]therefore we only need to show that $c_   \varepsilon\le c$ where $c$ does not depend on $\varepsilon,$ and this can be done arguing as in \cite[Section 6]{d2022multiple}. 
It remains only to show that the solution we found is not constant. Assume by contradiction that $u_0\in T_\varepsilon$ is a constant solution then \[
J_\varepsilon(u_0)=\frac{1}{4}\|u_0\|_{1,\varepsilon}^2+\left(\frac{1}{4}-\frac{1}{p}\right)|u_0|_{p,\varepsilon}^p=\frac{\mu_g(M)}{\varepsilon^3}\left[\frac{c_*^2}{4}+\left(\frac{1}{4}-\frac{1}{p}\right) c_*^p\right] \rightarrow+\infty \text { as } \varepsilon \rightarrow 0
\]that is a contradiction since \[
\max_{T_\varepsilon}J_\varepsilon\le C.
\]
\bibliographystyle{abbrv}
\bibliography{biblio.bib}
\end{document}